\documentclass[12pt,a4paper,english]{amsart}
\usepackage{babel}
\usepackage[latin1]{inputenc}
\usepackage{amsmath}
\usepackage{amssymb}

\theoremstyle{plain}
\newtheorem{Thm}{Theorem}[section]

\newtheorem{Lemma}[Thm]{Lemma}

\theoremstyle{definition}

\newtheorem{Def*}{Definition}

\newtheorem{Example}[Thm]{Example}
\newtheorem{Cor}[Thm]{Corollary}

\newcommand{\m}[1]{\mathrm{E}[#1]}
\newcommand{\rank}{\mathrm{rank}}
\newcommand{\cost}{\mathrm{cost}}
\newcommand{\exit}{\mathrm{exit}}
\newcommand{\comp}[1]{\overline{#1}}
\renewcommand{\O}{\emptyset}

\begin{document}

\title [The Buck-Chan-Robbins conjecture]
{A proof of a conjecture of Buck, Chan and Robbins on the random
assignment problem}

\author{Svante Linusson}
\address{Svante Linusson\\Department of Mathematics\\
   Link{\"o}pings universitet \\
   SE-581 83 Link{\"o}ping, Sweden}
\email{linusson@mai.liu.se}

\author{Johan W{\"a}stlund}
\address{Johan W{\"a}stlund\\ Department of Mathematics \\
   Link{\"o}pings universitet \\
   SE-581 83 Link{\"o}ping, Sweden}
\email{jowas@mai.liu.se}

\date{\today}

\begin{abstract}
We prove the main conjecture of the paper ``On the expected value of the
minimum assignment'' by Marshall~W.~Buck, Clara~S.~Chan, and
David~P.~Robbins (Random Structures \& Algorithms 21 (2002), no. 1, 33--58).
This is a vast generalization of a formula conjectured by Giorgio Parisi for
the $n$ by $n$ random assignment problem.
\end{abstract}

\maketitle

\section{The Parisi formula}
This work is motivated by a conjecture made in 1998 by the physicist Giorgio
Parisi \cite{P98}. Consider an $n$ by $n$ matrix of independent exp(1)
random variables. Parisi conjectured that the expected value of the minimal
sum of $n$ matrix entries of which no two belong to the same row or column,
is given by the formula \begin{equation} \label{Parisi}
\sum_{i=1}^n{\frac{1}{i^2}}.\end{equation} An equivalent setting is obtained
by considering the expected minimum cost of a perfect matching in a complete
$n$ by $n$ bipartite graph with independent $\exp(1)$ edge costs.

The problem had already received quite some attention. At the time, the main
open question was the value of the limit of the expected optimal value as
$n$ tends to infinity. A non-rigorous argument due to Marc~M\'ezard and
Parisi \cite{MP85} showed that the limit ought to be $\zeta(2)=\pi^2/6$.
David~Aldous \cite{A92} proved, using an infinite model, that the limit
exists. The striking conjecture of Parisi (obviously consistent with the
conjectured $\zeta(2)$-limit) paved the way for an entirely new approach. It
seemed likely that \eqref{Parisi} would yield to an inductive argument, and
that therefore the $\zeta(2)$-limit could be established by exact analysis
of the ``finite $n$'' case. The Parisi formula was almost immediately
generalized by Don Coppersmith and Gregory B. Sorkin \cite{CS98} to
\begin{equation}\label{Coppersmith-Sorkin} \sum_{\substack{i,j\geq
0\\i+j<k}}{\frac{1}{(m-i)(n-j)}}\end{equation} for $k$-assignments in an $m$
by $n$ matrix, $k\leq\min(m,n)$. It is not hard to verify that
\eqref{Coppersmith-Sorkin} specializes to \eqref{Parisi} when $k=m=n$.

The Coppersmith-Sorkin conjecture was then generalized in two different
directions by Marshall~W. Buck, Clara~S. Chan, and David~P. Robbins
\cite{BCR02} and by the present authors \cite{LW00}. In a recent paper
\cite{LW03}, we prove our own conjecture, thereby establishing the Parisi
and Coppersmith-Sorkin formulas. Remarkably, a proof of these formulas was
announced simultaneously by Chandra~Nair, Balaji~Prabhakar and Mayank~Sharma
\cite{NPS03}.

In this paper we state and prove a simultaneous generalization of the
Buck-Chan-Robbins conjecture and the main theorem of \cite{LW03}. Not only
is this a stronger result than that of \cite{LW03}, but it also provides a
considerable improvement of the proof. By combining the approach of
\cite{BCR02} with that of \cite{LW03}, we obtain a proof of Parisi's formula
which is shorter, simpler, and gives a far better insight into the problem
than our original proof.

\section{The Buck-Chan-Robbins formula}
In this section, we state the formula conjectured by Buck, Chan and Robbins
in \cite{BCR02}. This is a ``combinatorial'' formula, involving a binomial
coefficient. As is shown in \cite{BCR02}, there is an equivalent
``probabilistic'' version of the formula, and one can pass between the two
via elementary properties of the M\"obius function. In this paper we work in
the probabilistic setting. In order to prove our main theorem, Theorem
\ref{T:main}, it is therefore not necessary to take the detour through the
combinatorial formulas. The purpose of this section is therefore mainly to
introduce some notation and provide the background. We remark, however, that
both in \cite{LW00} and \cite{BCR02}, the discovery of the probabilistic
formulas was made through formal manipulation of the combinatorial formulas.
Therefore the latter have played an important role in obtaining the results
of this paper.

Let $M$ be an $m$ by $n$ matrix of nonnegative real numbers. The rows and
columns of the matrix will be indexed by \emph{weighted sets} $R$ and $C$
respectively. We may take $R=\{1,\dots,m\}$ and $C=\{1,\dots,n\}$, but the
sets $R$ and $C$ come with weight functions $w_R$ and $w_C$ respectively
that associate a positive weight to each element of the set. A
\emph{$k$-assignment} is a set $\pi\subseteq R\times C$ of $k$ matrix
positions, or \emph{sites}, of which no two belong to the same row or
column. An assignment will also be called an \emph{independent} set. The
\emph{cost} of $\pi$ is the sum $$\cost_M(\pi)=\sum_{(i,j)\in\pi}{M(i,j)}$$
of the matrix entries in $\pi$. We let ${\min}_k(M)$ denote the minimum cost
of all $k$-assignments in $M$.

In \cite{BCR02}, the Parisi and Coppersmith-Sorkin conjectures are
generalized to a certain type of matrix with entries which are exponential
random variables, but not necessarily with parameter 1. We say that a random
variable $x$ is \emph{exponential of rate $\alpha$} if $Pr(x>t) = e^{-\alpha
t}$ for every $t\geq 0$. In this case we write $x\sim \exp(\alpha)$. Buck,
Chan and Robbins considered the following type of matrix: For every
$(i,j)\in R\times C$, $M(i,j)$ is $\exp(w_R(i)w_C(j))$-distributed, and
independent of all other matrix entries. To state the formula, we use the
following notation: If $X$ is a set of rows, we let $w_R(X) = \sum_{i\in
X}{w_R(i)}$, and $w_R(\comp X) = w_R(R\backslash X)$. We use similar
notation for sets of columns.

\begin{Thm} [Conjectured by Buck, Chan and Robbins 2000] \label{T:BCR}
Let $M$ be a matrix as described above. Then
\begin{equation} \label{bcr}
\m{{\min}_k(M)} = \sum_{\substack{X\subseteq R\\ Y\subseteq C}}
{\binom{m+n-1-|X|-|Y|}{k-1-|X|-|Y|}
\frac{ (-1)^{k-1-|X|-|Y|}}{w_R(\comp X)
w_C(\comp Y)}}.
\end{equation}
\end{Thm}

Notice that in order for the binomial coefficient to be nonzero, we must
have $|X|+|Y|<k$, which ressembles the condition $i+j<k$ in the
Coppersmith-Sorkin formula \eqref{Coppersmith-Sorkin}. It is still not
entirely obvious that \eqref{bcr} specializes to the Coppersmith-Sorkin
formula when the row- and column weights are set to 1. However, in
\cite{BCR02}, the formula \eqref{bcr} is shown to be equivalent to a formula
given by an urn model. We will generalize this to a setting where a certain
set of matrix entries are set to zero.

\section{Main Theorem} \label{S:main}

The main theorem of \cite{LW03} is a formula for the expected value of the
minimal $k$-assignment in a matrix where a specified set of entries are set
to zero, and the remaining entries are independent $\exp(1)$-variables. In
this article we prove a formula for the common generalization of the
matrices considered in \cite{BCR02} and in \cite{LW03}. We say that $M$ is a
\emph{standard matrix} if the entries in a certain set $Z$ of sites are
zero, and the remaining entries are independent and distributed according to
the row- and column weights, that is, $M(i,j)\sim \exp(w_R(i)w_C(j))$. This
is an obvious generalization of the concept of standard matrix in
\cite{LW03}. As in \cite{LW00, LW03, BCR02} we give two seemingly different
but equivalent formulations of our main theorem.

Let $Z\subseteq R\times C$ be a set of sites. A \emph{file} is a row or a
column. Let $\lambda$ be a set of files. We say that $\lambda$ is a
\emph{cover} of $Z$ if every site in $Z$ lies in a file that belongs to
$\lambda$. By a cover of the matrix $M$ we mean a cover of the set of zeros
of $M$. By a \emph{$k-1$-cover} we mean a cover consisting of $k-1$ files.
Finally by a \emph{partial $k-1$-cover} we mean a subset of a $k-1$-cover.

Let $J_k(M)$ be the set of partial $k-1$-covers of the zeros of $M$. Let
$\hat J_k(M)$ denote the poset  consisting of $J_k(M)$ ordered by inclusion,
together with an artificial largest element $\hat 1$. The $k-1$-covers are
coatoms in $\hat J_k(M)$.
Let $\mu$ denote the M\"obius function on intervals in $\hat J_k(M)$ (see
e.g. \cite{S} for the basics of M\"obius functions).

The following is a combinatorial formulation of our main theorem.

\begin{Thm} [Main Theorem, combinatorial version] \label{T:main1}
Let $M$ be a standard matrix. Then
\begin{equation}\label{main1}
\m{{\min}_k(M)} = \sum_{(X,Y)\in J_k(M)}
\frac{-\mu((X,Y),\hat 1)}{w_R(\comp X)
w_C(\comp Y)}.
\end{equation}
\end{Thm}

If there are no zero entries in $M$, $J_k(M)$ consists of all sets of at
most $k-1$ files. The poset $\hat J_k(M)$ is a truncated Boolean lattice
obtained by deleting all elements of rank $\geq k$ except the top element in
the Boolean lattice $B_{m+n}$. The fact that \eqref{main1} specializes to
\eqref{bcr} follows from the fact that the M\"obius function of the
truncated Boolean lattice occurring in \eqref{main1} is given by the signed
binomial coefficient in \eqref{bcr}.

\begin{Example}
Let $M$ be a standard $2\times 2$ matrix with no zeros, and let the row- and
column weights be $w_R(i)=a_i$, and $w_C(j) = b_j$. With $k=2$, the poset
$\hat J_2(M)$ consists of six elements: The bottom element is the empty set.
There are four elements of rank 1 consisting of one file, and then there is
the top element $\hat 1$. The M\"obius function on the interval $(\O, \hat
1)$ is equal to 3, and the M\"obius function on the intervals from the rank
1 elements to the top element is $-1$.

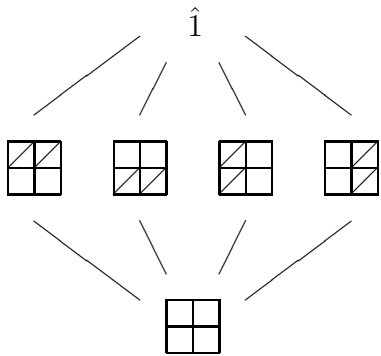
\begin{center}

\begin{picture} (160, 160) (00,0)
\put(70,10){\line(0, 1){20}}
\put(70,10){\line(1, 0){20}}
\put(80,10){\line(0, 1){20}}
\put(90,10){\line(0, 1){20}}
\put(70,20){\line(1, 0){20}}
\put(70,30){\line(1, 0){20}}

\put(10,70){\line(0, 1){20}}
\put(10,70){\line(1, 0){20}}
\put(20,70){\line(0, 1){20}}
\put(30,70){\line(0, 1){20}}
\put(10,80){\line(1, 0){20}}
\put(10,90){\line(1, 0){20}}

\put(50,70){\line(0, 1){20}}
\put(50,70){\line(1, 0){20}}
\put(60,70){\line(0, 1){20}}
\put(70,70){\line(0, 1){20}}
\put(50,80){\line(1, 0){20}}
\put(50,90){\line(1, 0){20}}

\put(90,70){\line(0, 1){20}}
\put(90,70){\line(1, 0){20}}
\put(100,70){\line(0, 1){20}}
\put(110,70){\line(0, 1){20}}
\put(90,80){\line(1, 0){20}}
\put(90,90){\line(1, 0){20}}

\put(130,70){\line(0, 1){20}}
\put(130,70){\line(1, 0){20}}
\put(140,70){\line(0, 1){20}}
\put(150,70){\line(0, 1){20}}
\put(130,80){\line(1, 0){20}}
\put(130,90){\line(1, 0){20}}

\put(10,80){\line(1, 1){10}}
\put(20,80){\line(1, 1){10}}
\put(50,70){\line(1, 1){10}}
\put(60,70){\line(1, 1){10}}
\put(90,70){\line(1, 1){10}}
\put(90,80){\line(1, 1){10}}
\put(140,70){\line(1, 1){10}}
\put(140,80){\line(1, 1){10}}

\put(60,30){\line(-4, 3){40}}
\put(100,30){\line(4, 3){40}}
\put(20,100){\line(4, 3){40}}
\put(140,100){\line(-4, 3){40}}
\put(70,40){\line(-1, 2){10}}
\put(90,40){\line(1, 2){10}}
\put(60,100){\line(1, 2){10}}
\put(100,100){\line(-1, 2){10}}

\put(78,130){\text{$\hat 1$}}
\end{picture}

\end{center}
Hence according to \eqref{main1}
\begin{multline} \label{example}
\m{{\min}_2(M)} =\\ \frac{-3}{(a_1+a_2)(b_1+b_2)} + \frac{1}{a_2(b_1+b_2)} +
\frac{1}{a_1(b_1+b_2)} + \frac{1}{b_2(a_1+a_2)} + \frac{1}{b_1(a_1+a_2)}.
\end{multline}
If we set the weights equal to 1, \eqref{example} specializes to $5/4$, in
accordance with the Parisi formula $1+1/4$. If we compare to \eqref{bcr}, we
see that the numerators in \eqref{example} are indeed equal to the binomial
coefficients in \eqref{bcr}
\end{Example}

\section{Matrix reduction}

Polynomial time algorithms for computing ${\min}_k(M)$ for a given
(non-random) matrix $M$ are well-known. We do not focus here on issues of
computational efficiency, but we outline an algorithm whose special features
will be of importance. The following lemma expresses a fundamental property
of optimal assignments. It is proved in \cite{BCR02}, although these authors
make no claims of originality. The first statement is certainly well-known,
but we haven't been able to trace the second statement to any other source
than \cite{BCR02}. For completeness, we include an outline of the proof.

\begin{Lemma}[Nesting Lemma] \label{L:nesting}
Let $M$be a real $m$ by $n$ matrix, and let $k_1\leq k_2\leq \min(m,n)$ be
positive integers. If $\mu$ is an optimal $k_1$-assignment in $M$, then
there is an optimal $k_2$-assignment $\mu'$ in $M$ such that every file that
intersects $\mu$ also intersects $\mu'$. Moreover, if $\nu$ is an optimal
$k_2$-assignment, then there is an optimal $k_1$-assignment $\nu'$ such that
every file that intersects $\nu'$ intersects $\nu$.
\end{Lemma}

\begin{proof}[Sketch of proof]
We may assume that $k_1=k_2-1$. Suppose that $\mu$ is an optimal
$k_1$-assignment, and $\nu$ is an optimal $k_2$-assignment. Consider the
symmetric difference $\delta = \mu\triangle \nu$. We say that two sites are
adjacent if they are in the same row or column. The components of $\delta$
with respect to adjacency are cycles or paths. Suppose that $\delta'$ is a
subset of $\delta$ which consists of a number of entire components of
$\delta$, in other words such that no site in $\delta'$ is adjacent to a
site in $\delta\backslash \delta'$. Suppose moreover that $\delta'$ is
balanced in the sense that it contains equally many sites from $\mu$ and
$\nu$. Then $\mu\triangle\delta'$ is a $k_1$-assignment, and
$\nu\triangle\delta'$ is a $k_2$-assignment. It follows that
$\cost(\mu\cap\delta') = \cost(\nu\cap\delta')$. Hence in the components of
$\delta$ that are balanced, the cost of the sites in $\mu$ is equal to the
cost of the sites in $\nu$, and for all the other components, the difference
in cost between the sites in $\mu$ and the sites in $\nu$ is the same (with
a sign depending on which one of $\mu$ and $\nu$ is overrepresented). As a
consequence, we may take a single component $\delta_1$ of $\delta$ so that
$\delta_1$ has one more site in $\nu$ than in $\mu$. Then $\mu' =
\mu\triangle\delta_1$ is an optimal $k_2$-assignment, and $\nu' =
\nu\triangle\delta_1$ is an optimal $k_1$-assignment, and it is
straightforward to verify that $\mu'$ and $\nu'$ have the desired
properties.
\end{proof}

Let $Z\subseteq R\times C$ be a set of sites. We say that a cover of $Z$ is
\emph{optimal} if it has the minimum number of files among all covers of
$Z$. The \emph{rank} of a set of sites is the size of the largest
independent subset. The following is a famous theorem due to Denes K\"onig:
\begin{Thm}[K\"onig's theorem]
The number of files in an optimal cover of $Z$ is equal to $\rank(Z)$.
\end{Thm}

The following theorem forms the basis of an algorithm for computing
${\min}_k(M)$:

\begin{Thm} \label{T:basis}
Let $M$ be a nonnegative $m$ by $n$ matrix, and let $\lambda$ be an optimal
cover of $M$. Suppose that there is no zero cost $k+1$-assignment in $M$.
Then every file in $\lambda$ intersects every optimal $k$-assignment in $M$.
\end{Thm}

\begin{proof}
Let $\mu$ be an optimal $k$-assignment. Since there is no zero cost
$k+1$-assignment, $k\geq |\lambda|$. By Lemma \ref{L:nesting}, there is an
optimal $|\lambda|$-assignment $\mu'$ such that $\mu$ intersects every file
that intersects $\mu'$. By K\"onig's theorem, there is a zero cost
$|\lambda|$-assignment, and since $\mu'$ is optimal, this means that $\mu'$
has zero cost. Hence every file in $\lambda$ intersects $\mu'$. The
statement follows.
\end{proof}

The following matrix operation is fundamental for the algorithm. We refer to
it as \emph{matrix reduction}. Let $M$ be a nonnegative $m$ by $n$ matrix,
and let $\lambda=X\cup Y$ be an optimal cover of $M$, where $X$ is the set
of rows and $Y$ is the set of columns in $\lambda$. The \emph{reduction}
$M'$ of $M$ by $\lambda$ is obtained from $M$ as follows: Let $t$ be the
minimum matrix entry of $M$ which is not covered by $\lambda$. If the site
$(i,j)$ is not covered by $\lambda$, we let $M'(i,j)=M(i,j)-t$. In
particular, this means that $M'$ will have a zero entry not covered by
$\lambda$. If the site $(i,j)$ is \emph{doubly covered} by $\lambda$, that
is, $i\in X$ and $j\in Y$, then we let $M'(i,j)=M(i,j)+t$. Finally if
$(i,j)$ is covered by exactly one file in $\lambda$, we let
$M'(i,j)=M(i,j)$. Notice that the entries of $M'$ are nonnegative.

\begin{Lemma} Let $M'$ be the reduction of $M$ by the optimal cover
$\lambda$. A $k$-assignment which is optimal in $M$ is also optimal in $M'$.
\end{Lemma}

\begin{proof}
Let $t$ be the minimum of the entries in $M$ that are not covered by
$\lambda$. For $s<t$, let $M_s$ be the matrix obtained from $M$ by
subtracting $s$ from the non-covered entries and adding $s$ to the doubly
covered entries. Since $M_s$ has no zero entries except those of $M$, it
follows from Theorem \ref{T:basis} that every optimal $k$-assignment in
$M_s$ must intersect every file of $\lambda$. By continuity, it follows that
there is some optimal $k$-assignment in $M'$ that intersects every file in
$\lambda$. All $k$-assignments that intersect every file of $\lambda$ are
affected in the same way by the reduction from $M$ to $M'$, namely if $\mu$
is such a $k$-assignment, then $\cost_\mu(M') =
\cost_\mu(M)-(k-|\lambda|)t$. Hence if $\mu$ is an optimal $k$-assignment in
$M$, then $\mu$ is optimal also in $M'$.
\end{proof}

{}From Lemma \ref{L:nesting} and K\"onig's theorem we can deduce the
following:

\begin{Lemma}
There is an optimal cover of $Z$ containing every row that belongs to some
optimal cover of $Z$, and similarly there is an optimal cover that contains
every column that belongs to some optimal cover.
\end{Lemma}

These covers are called the \emph{row-maximal} and the \emph{column-maximal}
optimal covers, respectively.

\begin{proof}
It follows immediately from K\"onig's theorem that a file belongs to an
optimal cover of $Z$ if and only if it intersects every maximal independent
subset of $Z$. Let $\lambda$ be the set of rows that belong to some optimal
cover of $Z$. Let $Z\backslash \lambda$ be the set of sites in $Z$ that are
not covered by $\lambda$. Let $\mu$ be a maximal independent subset of
$Z\backslash \lambda$. Then by Lemma \ref{L:nesting} there is a maximal
independent subset $\mu'$ of $Z$ which intersects every row that intersects
$\mu$. At the same time, $\mu'$ must intersect every row in $\lambda$.
Therefore $\rank(Z\backslash \lambda) = \rank(\mu) = \rank(Z)-|\lambda|$.
Hence $\lambda$ can be extended to an optimal cover of $Z$.
\end{proof}
We want to be able to do induction over matrix reduction. Therefore we need
the following lemma:

\begin{Lemma}
Let $M = M_0$ be a nonnegative $m$ by $n$ matrix, and let $k\leq \min(m,n)$.
For $i\geq 0$, let $M_{i+1}$ be the reduction of $M_i$ by the column-maximal
optimal cover of $M_i$. Then one of the matrices $M_i$ has a zero cost
$k$-assignment.
\end{Lemma}

\begin{proof}
Let $Z_i$ be the set of sites where $M_i$ has zeros. Let $\lambda_i$ be the
column-maximal optimal cover of $Z_i$. By K\"onig's theorem, $Z_0$ has an
independent subset $\mu$ containing exactly one site in each file in
$\lambda_0$. Hence $\mu$ contains no site which is doubly covered by
$\lambda_0$. It follows that $\rank(Z_1)\geq \rank(Z_0)$. Suppose that
$\rank(Z_1)=\rank(Z_0)$. Every column in $\lambda_1$ must belong to
$\lambda_0$. Consequently every row in $\lambda_0$ must belong to
$\lambda_1$. Now since $M_1$ has a zero which is not covered by $\lambda_0$,
there has to be a row in $\lambda_1$ which is not in $\lambda_0$.

To sum up, in each step of the reduction process, either the rank of the set
of zeros increases, or the number of rows in the column-maximal optimal
cover increases.
\end{proof}

A feature of matrix reduction that has been exploited in several papers
\cite{LW00, AS02} is that it keeps track of the cost of the optimal
$k$-assignment. In fact, if $t$ is as above and $M$ reduces to $M'$, then
${\min}_k(M) = (k-\left|\lambda\right|)\cdot t + {\min}_k(M')$. This means
that we can compute ${\min}_k(M)$ recursively by iterating the reduction and
keeping track of the values of $t$ as well as the sizes of the optimal
covers that are used. As long as the matrix entries are independent
exponential variables, it is easy to compute the expected value of the
minimum $t$, even for general $m$ and $n$. However, since the doubly covered
entries will eventually consist of sums of several dependent random
variables, it becomes extremely hard to reach any conclusions valid for
general $k$ through this approach.

One of the key insights that led to the proof of the Parisi formula in
\cite{LW03} was the fact that information about the probability that a
certain matrix element participates in the optimal assignment will give
information about the expected minimum cost. However, a problem with the
reduction algorithm is that in general, it loses track of the location of
the optimal assignment.

\begin{Example} Here $k=2$, and after the final step, the matrix contains
two zero-cost 2-assignments, of which only one was optimal in the original
matrix.

\begin{equation}\notag
\begin{pmatrix} 1 & 2 & 3 \\ 3 & 4 & 3 \end{pmatrix}
\overset{\O}\longrightarrow \begin{pmatrix} 0 & 1 & 2 \\ 2 & 3 & 2
\end{pmatrix} \overset{\{{\rm column 1}\}}\longrightarrow
\begin{pmatrix} 0 & 0 & 1 \\ 2 & 2 & 1 \end{pmatrix}
\overset{\{{\rm row 1}\}}\longrightarrow
\begin{pmatrix} 0 & 0 & 1 \\ 1 & 1 & 0 \end{pmatrix}
\end{equation}

\end{Example}

The approach taken in this paper builds on an observation that has largely
been overlooked, even in \cite{LW03}, namely that when the column-maximal
optimal cover is used, matrix reduction keeps track of the set of rows that
intersect the optimal $k$-assignment.

\section{The participation probability lemma}

In this section we prove a slightly refined version of a lemma which first
occurred in \cite{LW00}. This lemma describes the probability that a certain
exponential variable participates in the optimal solution to a random
assignment problem.

\begin{Lemma}[\cite{LW00}] \label{L}
Let $M$ be a random matrix where a particular entry $M(i,j)\sim\exp(\alpha)$
is independent of the other matrix entries. Let $M'$ be as $M$ except that
$M'(i,j)=0$. Then the probability that $(i,j)$ belongs to the optimal
$k$-assignment in $M$ is
$$\alpha\cdot\left(\m{{\min}_k(M)}-\m{{\min}_k(M')}\right).$$
\end{Lemma}

\begin{proof} We condition on all entries in $M$ except $M(i,j)$. Let $M_t$
be the deterministic matrix obtained by also conditioning on $M(i,j)=t$. Let
$f(t) = \min_k(M_t)$. Then either $f$ is constant, or $f$ increases linearly
up to a certain point after which it is constant. The key observation is
that the site $(i,j)$ belongs to the optimal $k$-assignment in $M_t$ if and
only if $f'(t)=1$ (disregarding the possibility that $t$ is equal to the
point where $f$ is not differentiable). Therefore if $x\sim \exp(\alpha)$,
then the probability that $(i,j)$ belongs to the optimal $k$-assignment in
$M_x$ is equal to $\m{f'(x)}$. By partial integration we have
\begin{multline} \m{f'(x)} = \alpha\int_0^\infty{e^{-\alpha t}f'(t) dt} =\\
\alpha\int_0^\infty{d\left(e^{-\alpha t}f(t)\right)} +
\alpha^2\int_0^\infty{e^{-\alpha t}f(t)dt} \\= - \alpha f(0) + \alpha
\m{f(x)} =
\alpha\cdot\left(\m{{\min}_k(M)}-\m{{\min}_k(M')}\right).\end{multline}
\end{proof}

\section{The Buck-Chan-Robbins urn model}
The following urn model is described in \cite{BCR02}: An urn contains a set
of balls, each with a given positive weight. Balls are drawn one at a time
without replacement, and each time the probability of drawing a particular
ball is proportional to the weight of the ball. This simple model has
perhaps been studied before, but the connection to random assignment
problems is due to Buck, Chan and Robbins.

To each weighted set we can associate such an urn process. Here we take as
our weighted set the set $R$ of row indices (in order not to make any secret
of the kind of application we have in mind). We consider a continuous time
version of this process. Each ball (row) $i$ remains in the urn for an
amount of time which is $\exp(w_R(i))$-distributed, and the times at which
the balls leave the urn are all independent.

The urn process is described by a continuous time random walk
$u_R:\mathbf{R}^+ \to 2^{R}$ on the power set of $R$. For $t\geq 0$,
$u_R(t)$ is the set of balls that have been drawn at time $t$.

If $X\subseteq R$, we denote by $Pr_R(X)$ the probability that this random
walk reaches $X$, that is, the probability that every ball in $X$ is drawn
before every ball not in $X$.
\begin{Example}
If there are three balls labeled $1, 2, 3$ then $$Pr_R(\O) = Pr_R(\{1, 2,
3\}) = 1,$$ $$Pr_R(\{1\}) = \frac{w_R(1)}{w_R(\{1, 2, 3\})}$$ and
$$ Pr_R(\{1, 2\}) = \frac{w_R(1)w_R(2)}{w_R(\{1, 2, 3\})w_R(\{2, 3\})}+
\frac{w_R(1)w_R(2)} {w_R(\{1, 2, 3\})w(\{1, 3\})},$$ since the set $\{1,
2\}$ can be obtained either by first choosing $1$ and then $2$, or the other
way around.
\end{Example}

By an \emph{order ideal} (or just ideal for short) we mean a family of sets
of balls which is closed under taking subsets.

\begin{Lemma}\label{L:leaveideal}
Suppose that $I$ is an order ideal and further that $\O\in I$ and $R\notin
I$. Then \begin{equation} \label{e789} \sum_{X\in
I}{\sum_{\substack{i\\X\cup \{i\} \notin I}}{\frac{w_R(i) Pr_R(X)}{w_R(\comp
X)}}} = 1.\end{equation}
\end{Lemma}

\begin{proof}
The random walk $u_R$ starts in $\O$ which is in $I$, and ends in $R$ which
is not in $I$. Since $I$ is an order ideal, there will be exactly one step
of the walk which leads from a set in $I$ to a set not in $I$. The left hand
side of \eqref{e789} sums the probabilities of leaving $I$ via a certain
step, taken over all possible ways of leaving $I$.
\end{proof}

If $I$ is an ideal, and $R\notin I$, then we let $\exit_R(I)$ be the random
subset of $R$ which is the first set in the urn process which does not
belong to $I$. For $i\in R$, we denote by $I\backslash i$ the ideal
consisting of all sets in $I$ which do not contain $i$.

\begin{Lemma} \label{L:exit} If $i\in R$ then $$Pr(i\in \exit_R(I)) =
\sum_{X\in I\backslash i}{\frac{w_R(i)Pr_R(X)}{w_R(\comp X)}}.$$
\end{Lemma}
\begin{proof}
$i\in \exit_R(I)$ if and only if the ball $i$ is drawn at a moment where the
set of balls already drawn belongs to $I$. Therefore we get the probability
of this event by summing the probability of first arriving at $X$ and then
drawing $i$ in the next step, over all sets $X$ in $I\backslash i$.
\end{proof}

We can also describe this probability inductively in terms of the
corresponding probability for smaller ideals. If $I$ is an ideal and $i\in
R$, we let $I/i = \{X\subseteq R: X\cup\{i\}\in I\}$.

\begin{Lemma} \label{L:exit2}
If $I$ is an ideal such that $\O\in I$ and $R\notin I$, and $i_0\in R$, then
$$Pr(i_0\in \exit_R(I)) =
\frac{w_R(i_0)}{w_R(R)}+\frac{1}{w_R(R)}\sum_{i\neq
i_0}{w_R(i)Pr(i_0\in\exit_R(I/i))}.$$
\end{Lemma}

\begin{proof}
The first term in the right hand side is the probability that $i_0$ is the
first ball to be drawn. The probability that ball $i$ is the first ball to
be drawn is $w_R(i)/w_R(R)$, and given that this is the case, the ball $i_0$
belongs to $\exit_R(I)$ if and only if it belongs to $\exit_R(I/i)$.
\end{proof}

If $I$ is an ideal, we let $T_R(I) = \inf(t: u_R(t)\notin I)$ denote the
\emph{exit time} of $I$, in other words the time at which the random walk
$u_R$ leaves $I$, or equivalently the amount of time it spends in $I$. The
following formula for the expected value $\m{T_R(I)}$ of the exit time
follows from the observation that the amount of time spent in $I$ is equal
to the sum of the time spent at each $X\in I$.
\begin{Lemma} \label{L:exittime}
$$\m{T_R(I)} = \sum_{X\in I}{\frac{Pr_R(X)}{w_R(\comp X)}}.$$
\end{Lemma}

\begin{proof}
Given that the walk reaches $X$, the expexted amount of time until another
ball is drawn is equal to $$\frac{1}{w_R(\comp X)}.$$
\end{proof}

\section{A formula for the participation probability of a row}
In this section we obtain a connection between the random assignment problem
and the urn model by deriving a formula for the probability that a certain
row intersects an optimal $k$-assignment. The special case of matrices
without zero entries was proved in \cite{BCR02}. Another special case, that
of row- and column-weights equal to 1 (rate 1 exponential variables) was
proved in \cite{LW03} by a different method.

If $M$ is a nonnegative random matrix, we let $\rho_k(M)$ be the (random)
set of rows that intersect some optimal $k$-assignment in $M$. If $Z$ is a
set of sites, we let $r(Z)$ be the set of rows in the row-maximal optimal
cover of $Z$. Moreover, if $k$ is a positive integer, we let $I_k(Z)$ be the
ideal of all sets of rows which are partial $k-1$-covers of $Z$.

\begin{Lemma} \label{L:extend}
$X\in I_k(Z)$ iff $X\cup r(Z)\in I_k(Z)$.
\end{Lemma}

\begin{proof}
We prove this by induction on the size of $X$. The induction step is
equivalent to proving that the statement holds when $X$ consists of one row,
say $X=\{i\}$. If $i\in r(Z)$, then $X=X\cup r(Z)$, so the statement is
obvious. If $i\notin r(Z)$, then let $Z'$ be the set of sites in $Z$ which
are not in row $i$. By K\"onig's theorem, $\rank(Z') = \rank(Z)$. Hence an
optimal cover of $Z$, in particular the row-maximal one, is also an optimal
cover of $Z'$.
\end{proof}

\begin{Cor} \label{C:contraction} Let $Z$ be a set of sites, and let $(i,j)$
be a site such that $\rank(Z\cup \{(i,j)\}) = \rank(Z) + 1$. Then
$I_k(Z\cup\{(i,j)\}) = I_k(Z)/i$.
\end{Cor}

\begin{proof}
Since $i\in r(Z\cup \{(i,j)\})$, this follows from Lemma \ref{L:extend}.
\end{proof}

The following theorem establishes the connection between the urn process and
the random assignment problem. Thereby it forms the basis for our approach,
and in a sense it is the central theorem in the paper.

\begin{Thm}
Let $M$ be an $m$ by $n$ random matrix indexed by weighted sets $R$ and $C$.
Suppose that $M$ has the following properties: There is a specified set $Z$
of sites where the entries in $M$ are zero. The remaining entries in the
rows in $r(Z)$ are positive real numbers. For $i\in R$, $j\in C$, if
$i\notin r(Z)$ and $(i,j)\notin Z$, then $M(i,j)$ is
$\exp(w_R(i)w_C(j))$-distributed and independent of the other matrix
entries. Suppose that a certain row $i_0$ has no zeros. Then $$Pr(i_0\in
\rho_k(M)) = Pr(i_0 \in \exit_R(I_k(Z))).$$
\end{Thm}

In our applications of this theorem, we are always dealing with standard
matrices. However, to make the inductive proof go through, we must condition
on the values of the nonzero entries in the rows in $r(Z)$. For this reason
we let these entries be fixed numbers instead of random variables.

\begin{proof}
We prove this by induction. Let $\lambda$ be the column-maximal optimal
cover of $Z$. Let $M'$ be the reduction of $M$ by $\lambda$. Then there is
at least one new zero in $M'$, that is, a site $(i,j)$ which is not covered
by $\lambda$ and such that $M'(i,j)=0$. We let $Z' = \{(i,j):M'(i,j) = 0\}$.
We consider two cases.

(1) All new zeros are in rows that belong to $r(Z)$. In this case $\rank(Z')
= \rank(Z)$, and consequently $r(Z') = r(Z)$. It follows immediately from
Lemma \ref{L:extend} that $I_k(Z') = I_k(Z)$. Hence by induction, $Pr(i_0\in
\rho_k(M)) = Pr(i_0\in\exit_R(I_k(Z)))$.

(2) There is a new zero $M'(i,j)$ such that $i \notin r(Z)$. Since $M(i,j)$
has continuous distribution and is independent of all other matrix entries,
we may assume that $M'(i,j)$ is the only new zero in $M'$. Since the site
$(i,j)$ is not covered by any optimal cover of $Z$, we have $\rank(Z')=1 +
\rank(Z)$. Hence $i\in r(Z')=r(Z)\cup \{i\}$.

If $i=i_0$, then every optimal $k$-assignment in $M'$ must intersect row
$i_0$. Since every optimal $k$-assignment in $M$ is optimal in $M'$, every
optimal $k$-assignment in $M$ must intersect row $i_0$. If on the other hand
$i\neq i_0$, then with probability 1, row $i_0$ participates either in all
or in none of the optimal $k$-assignments in $M'$. If we condition on the
values of $M'$ in row $i$, then $M'$ satisfies the criteria of the theorem.
Hence by induction, $Pr(i_0\in\rho_k(M')) = Pr(i_0\in \exit_R(I_k(Z')))$. By
Corollary \ref{C:contraction} we have $I_k(Z') = I_k(Z)/i$. Therefore if we
condition only on being in case 2, then \begin{multline} Pr(i_0\in
\rho_k(M)) =\\ \frac{w_R(i_0)}{w_R(\comp{r(Z)})} +\frac{1}{w_R(\comp
{r(Z)})}\sum_{\substack{i\notin r(Z)\\ i\neq i_0}}{w_R(i)Pr(i_0\in
\exit_R(I_k(Z)/i))},\end{multline} which by Lemma \ref{L:exit2} is equal to
$Pr(i_0\in \exit_R(I_k(Z)))$.
\end{proof}

\section{The two-dimensional urn-process}
At this point we introduce a kind of product of two urn processes. We
consider two independent urn processes on the weighted sets $R$ and $C$
respectively. What we here call the \emph{two-dimensional urn-process} is
just a piece of notation that makes it easy to state the generalization of
the Buck-Chan Robbins formula. The weight function is multiplicative: If
$X\subseteq R$ and $Y\subseteq C$, then we let $$w_{R\times C}(X,Y) =
w_R(X)w_C(Y).$$ Time is two-dimensional, and we let $$u_{R\times C}(x, y) =
(u_R(x), u_C(y)).$$ We further let $$Pr_{R\times C}(X,Y) = Pr_R(X)Pr_C(Y).$$
Since the two one-dimensional processes are statistically independent,
$Pr_{R\times C}(X,Y)$ is equal to the probability that there exists a point
$(x, y)$ in the time plane such that $u_{R\times C}(x, y) = (X, Y)$.

Let $J$ be an order ideal in $2^{R}\times 2^{C}$. In analogy with the
one-dimensional exit time, we define the two-dimensional exit time
$T_{R\times C}(J)$ to be the amount of two-dimensional time spent in $J$,
that is, the area of the region given by $u_{R\times C}(x, y)\in J$. We have
\begin{equation} \label{2dtime} \m{T_{R\times C}(J)} = \sum_{(X,Y)\in
J}{\frac{Pr_{R\times C}(X,Y)}{w_{R\times C}(\comp X, \comp
Y)}}.\end{equation} As indicated in the figure, given that the random
process reaches $(X,Y)$, the expected amount of time spent there is
$1/w_{R\times C}(\comp X, \comp Y)$.

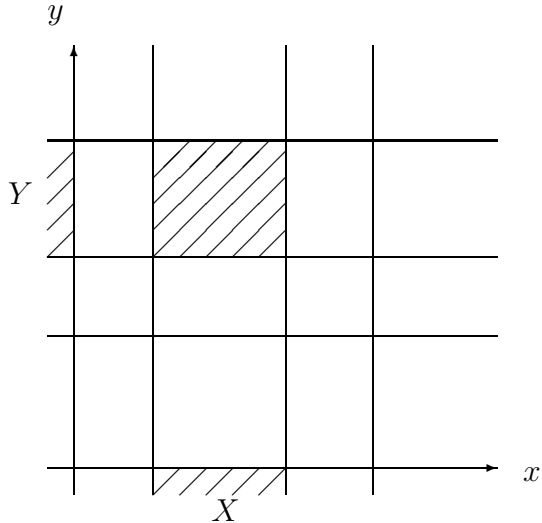
\setlength{\unitlength}{1.0pt}
\begin{center}
\begin{picture} (200, 250) (50,0)
\put(30,40){\vector(1, 0){170}} \put(40,30){\vector(0, 1){170}}
\put(70,30){\line(0, 1){170}}
\put(120,30){\line(0, 1){170}}
\put(153,30){\line(0, 1){170}}

\put(30,164){\line(1, 0){170}}

\put(30,90){\line(1, 0){170}}
\put(30,120){\line(1, 0){170}}

\put(70,120){\line(1, 1){44}}

\put(80,120){\line(1, 1){40}}
\put(90,120){\line(1, 1){30}}
\put(100,120){\line(1, 1){20}}
\put(110,120){\line(1, 1){10}}

\put(70,130){\line(1, 1){34}}
\put(70,140){\line(1, 1){24}}
\put(70,150){\line(1, 1){14}}

\put(70,30){\line(1, 1){10}}
\put(80,30){\line(1, 1){10}}
\put(90,30){\line(1, 1){10}}
\put(100,30){\line(1, 1){10}}
\put(110,30){\line(1, 1){10}}

\put(30,120){\line(1, 1){10}}
\put(30,130){\line(1, 1){10}}
\put(30,140){\line(1, 1){10}}
\put(30,150){\line(1, 1){10}}

\put(30,210){\text{$y$}}\put(210,35){\text{$x$}}
\put(15,140){\text{$Y$}}

\put(92,20){\text{$X$}}

\end{picture}
\end{center}
\setlength{\unitlength}{1pt}

\section{A formula for $\m{\min_k(M)}$}

Let $M$ be a standard matrix with rows and columns indexed by the weighted
sets $R$ and $C$. Let $J_k(M)$ be the ideal consisting of all partial
$k-1$-covers of the zeros of $M$. In this section we show that the expected
cost of the minimal $k$-assignment in $M$ is simply equal to the expected ex
it time of $J_k(M)$ in the two-dimensional urn process. When $M$ is a
matrix, we will write $I_k(M)$ for $I_k(Z)$, where $Z$ is the set of zeros
of $M$.

\begin{Thm}[Main Theorem, probabilistic version] \label{T:main}
\begin{equation}\label{cf} \m{{\min}_k(M)}=\m{T_{R\times C}(J_k(M))}.
\end{equation}
\end{Thm}

The proof of Theorem \ref{T:main} is inductive. We first prove that
\eqref{cf} is consistent with the row participation formula.

\begin{Lemma} \label{L:insertzero}
Let $M$ be a standard matrix where row $i_0$ contains no zeros. Let $M^{j}$
be obtained from $M$ by setting the entry in position $(i_0, j)$ equal to
zero. If $\m{{\min}_k(M^j)}=\m{T_{R\times C}(J_k(M^j))}$ for every $j$, then
$\m{{\min}_k(M)}=\m{T_{R\times C}(J_k(M))}$.
\end{Lemma}

\begin{proof}
By Lemma \ref{L}, the probability that the site $(i_0,j)$ belongs to the
optimal $k$-assignment in $M$ is
$$w_R(i_0)w_C(j)(\m{{\min}_k(M)}-\m{{\min}_k(M^j)}).$$ Therefore, summing
over $j$, $$w_R(i_0)\sum_{j\in
C}{w_C(j)\left(\m{{\min}_k(M)}-\m{{\min}_k(M^j)}\right)}=Pr(i_0\in\rho_k(M))
.$$

We divide by $w_R(i_0)$ and use the fact that by Lemma \ref{L:exit},
$$Pr(i_0\in\rho_k(M)) = Pr(i_0\in \exit_R(I_k(M))) = \sum_{X\in
I_k(M)\backslash i_0}{\frac{w_R(i_0)Pr_R(X)}{w(\comp X)}}.$$
Hence we obtain
\begin{equation} \label{eq:diff} \sum_{j\in C}
{w_C(j)\left(\m{{\min}_k(M)}-\m{{\min}_k(M^j)}\right)} = \sum_{X\in
I_k(M)\backslash i_0}{\frac{Pr_R(X)}{w(\comp X)}}.\end{equation}

It is clear that we can solve for $\m{{\min}_k(M)}$ in \eqref{eq:diff}. To
finish the induction step, it is therefore sufficient to prove that
\begin{equation} \label{toprove} \sum_{j\in C}{w_C(j)\left( \m{T_{R\times
C}(J_k(M))} - \m{T_{R\times C}(J_k(M^j))} \right)} = \sum_{X\in
I_k(M)\backslash i_0}{\frac{Pr_R(X)}{w_R(\comp X)}}.\end{equation}

If we fix a set $X\in I_k(M)\backslash i_0$, then by Lemma
\ref{L:leaveideal}, applied to the ideal $\{Y:(X,Y)\in J_k(M)\}$, we have

$$\sum_{\substack{Y\\ (X,Y)\in J_k(M)}}{\sum_{\substack{j\\(X,
Y\cup\{j\})\notin J_k(M)}}{\frac{w_C(j)Pr_C(Y)}{w_C(\comp Y)}}} = 1.$$

Since $(i_0,j)$ is the only zero of $M^j$ in row $i_0$, a set of files not
containing row $i_0$ can be extended to a $k-1$-cover of the zeros of $M^j$
if and only if this can be done while making use of column $j$. If we want
to cover the zeros of $M^j$ as efficiently as possible, there is no point in
using row $i_0$ if instead we can use column $j$. Therefore the condition
$(X, Y\cup\{j\})\notin J_k(M)$ on $j$ in the inner sum can be replaced by
$(X, Y)\notin J_k(M^j)$.

If we multiply by $Pr_R(X)/w_R(\comp X)$ and sum over all $X\in
I_k(M)\backslash i_0$, we see that the right hand side of \eqref{toprove}
equals

$$\sum_{X\in I_k(M)\backslash i_0}{\sum_{\substack{Y\\ (X,Y)\in
J_k(M)}}{\sum_{\substack{j\\(X, Y)\notin
J_k(M^j)}}{\frac{w_C(j)Pr_R(X)Pr_C(Y)}{w_R(\comp X)w_C(\comp Y)}}}}.$$
Here we can drop the conditions on $X$, since the inner sum will be empty
unless $X\in I_k(N)\backslash i_0$. After changing the order of summation so
that the sum over $j$ becomes the outer sum, this is equal to $$\sum_{j\in
C}{w_C(j)\sum_{\substack{(X,Y)\in J_k(M)\\(X,Y)\notin
J_k(M^j)}}{\frac{Pr_{R\times C}(X,Y)}{w_{R\times C}(\comp X, \comp Y)} }}
.$$

By \eqref{2dtime}, this is equal to the left hand side of \eqref{toprove}.
\end{proof}

Secondly, we show that \eqref{cf} is consistent with removing a column that
contains at least $k$ zeros.

\begin{Lemma} \label{L:removecolumn}
Suppose that $\m{{\min}_{k-1}(M)}=\m{T_{R\times C}(J_{k-1}(M))}$ for every
standard matrix $M$, in other words, suppose that \eqref{cf} holds when $k$
is replaced by $k-1$. Let $M$ be a standard matrix that has a column with at
least $k$ zeros. Then $\m{{\min}_k(M)}=\m{T_{R\times C}(J_k(M))}$.
\end{Lemma}

\begin{proof} Suppose that column $j_0$, has at least $k$ zeros. Let $M'$ be
the $m$ by $n-1$ matrix obtained by deleting the column $j_0$ of $M$. Since
every $k-1$-assignment in $M'$ can be extended to a $k$-assignment in $M$ by
including a zero in column $j_0$, we have
$\m{{\min}_k(M)}=\m{{\min}_{k-1}(M')}$. To prove the lemma, we therefore
show that with the obvious coupling of the urn processes, $T_{R\times
C}(J_k(M)) = T_{R\times C\backslash\{j_0\}}(J_{k-1}(M'))$.

Since there are $k$ zeros in column $j_0$, every $k-1$-cover of $M$ must
include $j_0$. Therefore $(X,Y)\in J_k(M)$ if and only if
$(X,Y\backslash\{j_0\})\in J_{k-1}(M')$. The lemma follows.
\end{proof}

We are now in a position to prove that \eqref{cf} holds whenever $m$ is
sufficiently large compared to $k$.
\begin{Thm}
If $M$ is a standard $m$ by $n$ matrix with $m>(k-1)^2$, then
$\m{{\min}_k(M)} = \m{T_{R\times C}(J_k(M))}$.
\end{Thm}

\begin{proof}
By Lemmas \ref{L:insertzero} and \ref{L:removecolumn}, it is sufficient to
prove that the statement holds when $M$ has at least one zero in each row,
and no column with $k$ or more zeros. In this case each column can contain
the leftmost zero of at most $k-1$ rows. Since there are more than $(k-1)^2$
rows, there must be at least $k$ columns that contain the leftmost zero of
some row. This implies that there is a zero cost $k$-assignment in $M$.
Consequently there is no $k-1$-cover, that is, $J_k(M) = \O$. Plainly
$\m{{\min}_k(M)} = 0 = T_{R\times C}(\O) = \m{T_{R\times C}(J_k(M))}$.
\end{proof}

Finally we prove that \eqref{cf} holds also for smaller matrices by taking
the limit as the weights of the exceeding rows tend to zero.

\begin{proof}[Proof of Theorem \ref{T:main}]
We prove \eqref{cf} by downwards induction on the number of rows. Suppose
that $M$ is a standard $m$ by $n$ matrix. Let $M_\epsilon$ be an augmented
matrix of $m+1$ rows and $n$ columns, so that the first $m$ rows equal $M$,
and row $m+1$ has no zeros and weight $w_R(m+1) = \epsilon$. When $\epsilon$
is small, the entries of row $m+1$ are large. With high probability, none of
them will participate in the optimal $k$-assignment, and consequently
$$\m{{\min}_k(M)}=\lim_{\epsilon\to 0}{\m{{\min}_k(M_\epsilon)}}.$$
We therefore have to show that $$\m{T_{R\times C}(J_k(M))} =
\lim_{\epsilon\to 0}{\m{T_{R\times C}(J_k(M_\epsilon))}}.$$

Since row $m+1$ does not have any zero, we have $T_{R\times
C}(J_k(M_\epsilon)) \leq T_{R\times C}(J_k(M))$ under the obvious coupling
of the corresponding urn processes. Hence we can squeeze $\m{T_{R\times
C}(J_k(M_\epsilon))}$ by $$(1-p)\m{T_{R\times C}(J_k(M))} \leq \m{T_{R\times
C}(J_k(M_\epsilon))} \leq \m{T_{R\times C}(J_k(M))},$$ where $p$ is the
probability that there is a point $(x, y)$ in the time plane such that
$u_R(x), u_C(y)$ is a partial $k-1$-cover for $M$ but not for $M_\epsilon$.
The only way we can have a partial $k-1$ cover for $M$ and not for
$M_\epsilon$ is if the row $m+1$ has been drawn. However, as $\epsilon\to
0$, the probability that row $m+1$ is drawn among the first $k-1$ rows goes
to zero. Hence as $\epsilon$ tends to zero, so does $p$. The theorem
follows.
\end{proof}

\section{Proof of the Buck-Chan-Robbins formula}

We prove that the two formulations of the main theorem are indeed
equivalent. We have \begin{multline} \sum_{(X,Y)\in J_k(M)}
\frac{-\mu((X,Y),\hat 1)}{w_R(\comp X)w_C(\comp Y)} =\\
\sum_{(X,Y)\in J_k(M)}{(-\mu((X,Y),\hat 1))\cdot \m{T_R(\{X'\subseteq
X\})}\cdot \m{T_C(\{Y'\subseteq Y\})}},\end{multline} since the factors
$$\frac{1}{w_R(\comp X)},\quad \frac{1}{w_C(\comp Y)}$$ can be interpreted
as the expected exit times of the ideals $\{X'\subseteq X\}$ and
$\{Y'\subseteq Y\}$ of all subsets of $X$ and $Y$ respectively. By Lemma
\ref{L:exittime},

\begin{multline}
\sum_{(X,Y)\in J_k(M)}
\frac{-\mu((X,Y),\hat 1)}{w_R(\comp X)w_C(\comp Y)}=\\
\sum_{(X,Y)\in J_k(M)}
-\mu((X,Y),\hat 1) \sum_{X'\subseteq X} \frac{Pr_R(X')}{w_R(\comp X')}
\sum_{Y'\subseteq Y} \frac{Pr_C(Y')}{w_C(\comp Y')}
\end{multline}
We now change the order of summation and get
\begin{equation} \label{2895}
\sum_{(X',Y')\in J_k(M)}
\frac{Pr_R(X')Pr_C(Y')}{w_R(\comp X')w_C(\comp Y')}
\sum_{(X',Y')\le (X,Y)< \hat 1} -\mu((X,Y),\hat 1).\end{equation}
In the factor to the right, we are summing over all $(X,Y)$ that satisfy
$(X',Y')\le (X,Y)< \hat 1$ in the poset $\hat J_k(M)$. By the definition of
the M\"obius function, this sum is equal to 1, so that we can drop this
factor. Hence \eqref{2895} is equal to
$$
\sum_{(X',Y')\in J_k(M)}
\frac{Pr_R(X')Pr_C(Y')}{w(\comp X')w(\comp Y')} = \m{T_{R\times C}(J_k(M))}.
$$

This specializes to the equivalence of the two formulations of the main
conjectures in \cite{LW00} and of the two formulations in \cite{BCR02}. The
more general setting here allows us to give a shorter proof.

To finish the proof of Theorem \ref{T:BCR}, we cite a well-known theorem
(G-C Rota?) that states that the M\"obius function of the truncated Boolean
lattice is indeed given by the binomial coefficient occurring in
\eqref{bcr}.

\begin{Thm}
Let $P$ be a poset consisting of the elements of rank $0, \dots, k-1$ in a
Boolean lattice of degree $N$, together with an artificial top element $\hat
1$. Then $$\mu(\O, \hat 1) = (-1)^k \binom{N-1}{k-1}.$$
\end{Thm}

\begin{proof}[Sketch of proof]
Every interval of the form $(\O, x)$ for $x\neq \hat 1$ is Boolean, and
therefore $\mu(\O, x) = (-1)^{\rank(x)}$. Hence $$\mu(\O, \hat 1)
= -\sum_{\O\leq x<\hat 1}{(-1)^{\rank(x)}}
= -\sum_{i=0}^{k-1}{(-1)^i\binom{N}{i}}.$$ The theorem now follows
immediately by induction on $k$.
\end{proof}

\section{Specializing to rate 1 variables}

In this section we briefly comment on the implications of Theorem
\ref{T:main} to the case of rate 1 variables. In particular, we show that
the Coppersmith-Sorkin formula follows. When there are no zeros in the
matrix, the ideal $J_k(M)$ consists of all sets of at most $k-1$ files.
Moreover, when all files have the same weight, the files become
indistinguishable in the urn process. The expected amount of time until the
next ball is drawn depends only on the number of balls already drawn, and
not on which particular balls have been drawn. Therefore we obtain the
expected value of $T(J_k(M))$ by simply conditioning on every step in the
urn process taking exactly its expected amount of time. The expected amount
of two-dimensional time that the process spends at the point where exactly
$i$ rows and $j$ columns have been drawn is equal to
$$\frac{1}{(m-i)(n-j)}.$$ {}From this, the Coppersmith-Sorkin formula follows.

As $n\to\infty$, the process can be approximated by continuous exponential
decay. At time $t$, a fraction of $e^{-t}$ of the balls will remain in the
urn. In the two-dimensional process, the borderline at which the process
exits the ideal $J_n(M)$ will approach the curve given by the equation
$$e^{-x} + e^{-y} = 1.$$ Hence we may obtain the limit value as the area
under this curve, which is indeed equal to $\zeta(2)$.

\end{document}